\documentstyle[12pt,twoside]{amsart}

\renewcommand{\c}[0]{{\mathbb C}}

\renewcommand{\r}[0]{{\mathbb R}} 

\renewcommand{\a}[0]{{\mathbb A}}

\newcommand{\p}[0]{{\mathbb P}}

\newcommand{\map}[0]{\dasharrow}
\newcommand{\qtq}[1]{\quad\mbox{#1}\quad}

\newcommand{\rank}[0]{\operatorname{rank}}
\newcommand{\mult}[0]{\operatorname{mult}}

\newcommand{\jac}[0]{\operatorname{Jac}}    
\newcommand{\dist}[0]{\operatorname{dist}}

\newsymbol\subsetneq 2328
\newsymbol\onto 1310

\newsymbol\twoheadrightarrow 1310

\newtheorem{thm}{Theorem}

\newtheorem{lem}[thm]{Lemma}
\newtheorem{cor}[thm]{Corollary}

\newtheorem{prop}[thm]{Proposition}

\theoremstyle{definition}
\newtheorem{defn}[thm]{Definition}

\newtheorem{say}[thm]{}
\newtheorem{exmp}[thm]{Example}

\newtheorem{rem}[thm]{Remark}          
\newtheorem{ack}{Acknowledgments}

\theoremstyle{remark}

\begin{document}
\bibliographystyle{amsplain}

\title{An Effective  \L ojasiewicz Inequality for\\
 Real Polynomials}
\author{J\'anos Koll\'ar}

\maketitle

Let $f_1,\dots,f_m$ be real analytic functions on $\r^n$
and $Z:=(f_1=\cdots=f_m=0)$ their common zero set.
The classical theorem of \cite{loj-stud} says that 
for every $P\in Z$ 
there are constants
$C,\epsilon, M>0$ such that
$$
\max_i\{|f_i({\bold x})|\}\geq C\cdot \dist({\bold x}, Z)^M\qtq{if
$\dist({\bold x}, P)<\epsilon$.}
$$
The smallest such $M$ is called the 
{\it \L ojasiewicz exponent} of $f_1,\dots,f_m$ at $P$. 
The computation or estimation of the
\L ojasiewicz exponent is a quite interesting problem.
For instance, if the $f_i$ are polynomials of degree $\leq d$,
one would like to have an explicit bound for $M$ in terms of $d$
and $n$. The following example shows the magnitude of the bounds
one can expect:

\begin{exmp}\label{worst.exmp}
  Set $f_1=x_1^d$ and $f_i=x_{i-1}-x_i^d$
for $i=2,\dots, n$. Then $\Phi({\bold x}):=\max_i\{|f_i({\bold
x})|\}>0$ for
${\bold x}\neq 0$.
Let $p(t)=(t^{d^{n-1}}, t^{d^{n-2}},\dots,t)$.
Then $\lim_{t\to 0}||p(t)||/|t|=1$ and $\Phi(p(t))=t^{d^n}$. 
Thus the \L ojasiewicz exponent is $\geq d^n$.
(In fact it equals $d^n$.) This works both over $\r$ and $\c$.

In the real case set $F=\sum f_i^2$. Then  $\deg F=2d$, $F$ has an
isolated real zero at the origin and the 
\L ojasiewicz exponent is $2d^n$.
\end{exmp}

The complex analytic variant of this question has been settled
in the papers \cite{null, jks, ckt}. It seems, however,
more difficult to obtain   effective estimates in the real
case.   \cite{solerno} works in full generality, but
the resulting bounds are quite large. 
For a single real polynomial \cite{gwoz1} established
the following  very nice estimate.
In view of (\ref{worst.exmp}), it is close to being optimal.

\begin{thm}\label{gwoz.loj.thm}\cite{gwoz1}
 Let $f(x_1,\dots,x_n)$ be a real polynomial of degree
$d$ with a strict local minimum at $f(0)=0$. Then there are
constants
$C,\epsilon>0$ such that
$$
f({\bold x})\geq C\cdot ||{\bold x}||^{(d-1)^n+1}\qtq{for
$||{\bold x}||<\epsilon$.}
$$
\end{thm}

A similar estimate for ${\bold x}$ near infinity is given
in \cite{gwoz2}.  These  give  estimates in case of several
functions
$f_i$ by considering the single function $\sum f_i^2$ instead.

The aim of this paper is to give an upper bound
for the \L ojasiewicz exponent  with
$\max_i\{|f_i({\bold x})|\}$  replaced by 
 $\max_i\{f_i({\bold x})\}$. 
This version is more general and
as an application we also get an
effective \L ojasiewicz estimate in case  ${\bold x}$ runs through
  a semi--algebraic set instead of the whole  $\r^n$.
Unfortunately, I have to assume that 
$\max_i\{f_i({\bold x})\}$ has an isolated zero.

\begin{defn} For any $n$ let $B(n)$ denote the largest
of the binomial coefficients $\binom{n}{k}$. That is
$$
B(n)=\binom{2m}{m} \qtq{if $n=2m$  and}
B(n)=\binom{2m+1}{m} \qtq{if $n=2m+1$.}
$$
\end{defn}

\begin{thm}[Effective \L ojasiewicz
inequality]{\ }\label{real.loj.max.thm}

 Let $f_i: i\in I$ be finitely many 
real polynomials of degree $\leq d$ in $n$ variables.
Set $\Phi({\bold x}):=\max_i\{f_i({\bold
x})\}$. 
\begin{enumerate}
\item Assume that $\Phi({\bold x})>0$ for 
$0<||{\bold x}||\ll 1$. Then there are constants $C,\epsilon$
such that
$$
\Phi({\bold x})\geq C\cdot ||{\bold x}||^{B(n-1)d^n}
\qtq{for $||{\bold x}||<\epsilon$.}
$$

\item Assume that $\Phi({\bold x})>0$ for 
$ ||{\bold x}||\gg 1$. Then there are constants $C,N$
such that
$$
\Phi({\bold x})\geq C\cdot ||{\bold x}||^{-B(n-1)d^n}
\qtq{for $||{\bold x}||>N$.}
$$
\end{enumerate}
\end{thm}

\begin{rem} It is possible that in (\ref{real.loj.max.thm})
the exponent $B(n-1)d^n$ can be replaced by $d^n$.
This would be   optimal   in view of
(\ref{worst.exmp}). In the case  of complex variables
the optimal exponent is $d^n$ by \cite[4.1]{null}.

If the polynomials $f_i$ have different degrees then
in the complex case one can use the B\'ezout--type bound
$\prod \deg f_i$ instead of $d^n$. In the real case
this does not hold any more. For instance,
set $f_1:=F$ as in (\ref{worst.exmp})
and $f_2=\cdots=f_{n+1}=x_{n+1}$.

It would be also of interest to find the best exponent 
in the case of a single polynomial as in (\ref{gwoz.loj.thm}).
The polynomial $F$ constructed in (\ref{worst.exmp})
is not always the worst one. For instance, the examples
of \cite{yoshi} give a degree $6$
polynomial in 2 variables whose \L ojasiewicz exponent
is $20>2\cdot 3^2$. It is quite likely that
K3 surfaces give similar examples with
$n=3, d=4$ and \L ojasiewicz exponent $20>2\cdot 2^3$.
Pemantle pointed out that if 
$F(x_1,\dots,x_n)$ has \L ojasiewicz exponent $L$
then $F(x_1,\dots,x_n)+(\ell(x_1,\dots,x_n)-x_{n+1}^d)^2$
has \L ojasiewicz exponent $dL$ for 
 a suitable linear form $\ell$. 
This gives an infinite sequence 
of degree 6 polynomials
whose \L ojasiewicz exponent is bigger than those in 
 (\ref{worst.exmp}).
I do not know any other examples
which are worse than (\ref{worst.exmp}).

In all the examples where the \L ojasiewicz exponent
is large, the function has a type $A$ critical point.
That is, in suitable analytic coordinates it is
of the form $z_1^2+\cdots+z_{n-1}^2+z_n^M$
and $M$ is the \L ojasiewicz exponent. 
(Such a singularity is said to be of type $A_{M-1}$.)
I do not know
if there are other types of extremal examples.
In any case, this leads to the following interesting problem:

Let $X\subset \p^n$ be a hypersurface of degree $d$
with an isolated $A_m$ singularity. Find an upper bound 
for $m$ in terms of $d,n$ and find examples where $m$ is large.
Upper bounds which are smaller than $d^n$
are given for surfaces by 
 \cite{miya} ($\frac23d^3$) and by \cite{var}
in general ($\frac{23}{24}d^n$).
\end{rem}

It is also of interest to study similar questions
when $\bold x$ runs through the points of a closed
semi--algebraic set $X$. 
For instance a classical result of P\'olya   (cf.\ \cite[2.24]{hlp})
studies the case when $X$ is the positive quadrant.
In these situations it is necessary to control
the degree of the defining inequalities of $X$.

\begin{cor}[Effective semi--algebraic \L
ojasiewicz inequality]\label{semialg.loj.cor}{\ }

Let $X$ be a closed  semi--algebraic set of the form 
$$
X:=\{{\bold x}\in  \r^n | g_j({\bold x})=0, h_k({\bold x})\geq
0\},
$$
 where the $g_j,h_k$ are real polynomials of degree   $\leq d$.
 Let $f_i$ be finitely many 
real polynomials of degree $\leq d$ in $n$ variables.
Set $\Phi({\bold x}):=\max_i\{f_i({\bold
x})\}$. 
\begin{enumerate}
\item Assume that $\Phi({\bold x})>0$ if ${\bold x}\in X$ and
$0<||{\bold x}||\ll 1$. Then there are constants $C,\epsilon$
such that
$$
\Phi({\bold x})\geq C\cdot ||{\bold x}||^{B(n-1)d^n}
\qtq{if ${\bold x}\in X$ and $||{\bold x}||<\epsilon$.}
$$

\item Assume that $\Phi({\bold x})>0$ if ${\bold x}\in X$ and
$ ||{\bold x}||\gg 1$. Then there are constants $C,N$
such that
$$
\Phi({\bold x})\geq C\cdot ||{\bold x}||^{-B(n-1)d^n}
\qtq{if ${\bold x}\in X$ and $||{\bold x}||>N$.}
$$
\end{enumerate}
\end{cor}

Proof. Set $\Psi({\bold x})=\max\{f_i({\bold x}), g_j({\bold x}),
-g_j({\bold x}), -h_k({\bold x})\}$. Then
$\Psi({\bold x})>0$ for every ${\bold x}\not\in X$ and
$\Psi({\bold x})=\Phi({\bold x})$ for ${\bold x}\in X$. 
Thus (\ref{real.loj.max.thm}) applies to $\Psi$
and we get (\ref{semialg.loj.cor}).\qed

\begin{say}[Proof of (\ref{real.loj.max.thm})]{\ }

Following a suggestion of Gwo\'zdziewicz, we use the
equivalent norm $\max_i|x_i|$ instead of $||{\bold x}||$.
The formula given by (\ref{crit.pt.formula}) becomes very explicit
(\ref{expl.formula}) and this
 leads to a better estimate in (\ref{real.loj.max.thm}).

Assume that we are in the situation of (\ref{real.loj.max.thm}.1).
Choose $r>0$ such that $0$ is the only zero of $\Phi$ 
in the cube 
$\max_i|x_i|\leq r$ and $\eta^*>0$ such that $\Phi({\bold x})\geq
\eta^*$   if $\max_i|x_i|= r$.

Consider the level set  $X_{\eta}:=\{{\bold x} | \Phi({\bold x})=\eta\}$.
We usually think of $X_{\eta}$ as a subset of $\c^n$.
The set of its real points is denoted by $ X_{\eta}(\r)$.
For $\eta^*>\eta> 0$ let 
$M_{\eta}\subset X_{\eta}(\r)$ 
denote the union of those  connected components  
which are contained in the cube 
$\max_i|x_i|\leq r$.
 If   $X_{\eta}$ is in ``general
position" then $M_{\eta}$ has a unique point $P_{\eta}\in
M_{\eta}$ furthest from the origin. The points
$P_{\eta}$ form a real curve $C_{\Phi}$ and it is sufficient to
prove that
$$
\Phi({\bold x})\geq C\cdot ||{\bold x}||^{B(n-1)d^n}
\qtq{for ${\bold x}\in C_{\Phi}$.}
$$
We expect that $C_{\Phi}$ is an open subset of the real points of
a complex algebraic curve $C^*$.  
If $D$ is any algebraic curve  then an easy argument
 shows that 
$$
\Phi({\bold x})\geq C\cdot ||{\bold x}||^{d\cdot \deg D}
\qtq{for ${\bold x}\in D$.}\qquad\qquad\qquad (*)
$$
The first part of the proof establishes the existence of the
curve $C^*$. This can not be done for every $\Phi$
but we see that $C^*$ exists after a suitable small perturbation
of $\Phi$. In the second part we  
use a slight strenghtening of the estimate ($*$)
to conclude the proof.

The proof of 
 (\ref{real.loj.max.thm}.2) is entirely analogous but instead of
working in a small cube we work in the complement of
a large cube.
\end{say}

\begin{say}[Defining $C^*$]{\ }
 For $\eta^*>\eta> 0$ let 
$h(\eta)$ be  the maximum of $\max_i|x_i|$ on $M_{\eta}$ and
$\Pi_{\eta}\subset
M_{\eta}$  the set of points where 
$\max_i|x_i|=h(\eta)$.
The union of these sets $\Pi=\cup_{\eta} \Pi_{\eta}\subset \r^n\times
(0,\eta^*)$ is  semi--algebraic.

Assume first that there are $i\neq j$ and a sign $\epsilon=\pm 1$ such
that
$\Pi_{\eta}\cap (x_i=\epsilon x_j)$ is nonempty for every $0<\eta\ll
\eta^*$. Then one can find $C^*$ inside the 
 hyperplane $(x_i=\epsilon x_j)\cong \r^{n-1}$.
This leads to   a much better bound.

Otherwise, up to permuting the coordinates, we may assume that
$\Pi_{\eta}$ contains an interior point of the face $x_1=h(\eta)$
of the cube $\max_i|x_i|=h(\eta)$ for every $0<\eta\ll
\eta^*$. 
 Let 
$$
J:=\{j| f_j({\bold x})=\eta\ \forall {\bold x}\in 
\Pi_{\eta}\cap (x_1=h(\eta))\},
$$
and set 
$$
X_{J,\eta}:=\{{\bold x}|f_j({\bold x})=\eta\ \forall j\in
J\}.
$$
 Then every point of $\Pi_{\eta}$ is a 
(not necessarily strict) local minimum 
 of $x_1$ restricted
to $X_{J,\eta}$.

Let $\jac_J=\jac(x_1, f_j:j\in J) $  denote the
Jacobian matrix of the functions $x_1$ and $f_j:j\in
J$.

 Note that $J$ does not depend on $\eta$.
If $X_{J,\eta}$ is smooth of codimension $|J|$ then
the critical points of $x_1$ on $X_{J,\eta}$
are defined by the system of equations
$$
f_j({\bold x})=\eta\ \forall j\in J\qtq{and} \rank (\jac_J)\leq
|J|.
$$
If we let $\eta$ vary then we obtain the family of critical points
defined by the equations
$$
f_j({\bold x})=f_{j'}({\bold x})\ \forall j,j'\in J\qtq{and} \rank
(\jac_J)\leq |J|.
$$
For sufficiently general $f_i$ this gives  exactly
the expected curve $C^*$ but if $X_{J,\eta}$ has
singular points or components of the wrong dimension then these
are counted as critical points by the above equations.
\end{say}

The next two elementary lemmas establish that these problems can
be avoided by a small perturbation of the $f_j$ and of the
cordinate system.

\begin{lem} Let $f_i: i\in I$ be finitely many 
real polynomials. Then there are constants $1/2<\epsilon_i<2$
and $0<\eta_0$ such that for every $J\subset I$ and
$0<|\eta|<\eta_0$ the variety
$$
X_{J,\epsilon_i,\eta}:=\{{\bold x}\in \c^n: \epsilon_jf_j({\bold
x})=\eta
\ \forall j\in J\}
$$
is a smooth complete intersection of codimension $|J|$ (or it is empty).
\end{lem}

For the proof it is quite important that 
$X_{J,\epsilon_i,\eta}$ is viewed as a  {\it complex} algebraic variety.
The result holds if we allow complex values of $\eta$
as well, but we do not need this.

Proof. Let $B\subset \c^{|I|+1}$ be the set of those
$(\epsilon_i,\eta)$ such that $X_{J,\epsilon_i,\eta}$
is either singular or has a component of 
codimension $<|J|$ for some $J$. $B$ is a constructible
 algebraic subset of $\c^{|I|+1}$  whose Zariski closure $\bar B$
is different from $\c^{|I|+1}$  by the Bertini theorem.
Choose a point $(\epsilon_i,\eta')\not\in \bar B$
satisfying $1/2<\epsilon_i<2$. 
Only finitely many points of the disc
$\{(\epsilon_i,t)|0<|t|<\eta'\}$ belong to $\bar B$.
Choose $\eta_0$ to be the smallest value which gives  a point in
$\bar B$.\qed

\medskip

The following lemma is   known in algebraic geometry
as the existence of Lefschetz pencils. A very readable self contained 
proof is given in the lectures \cite[1.8]{loo-pc}.

\begin{lem} Let $X\subset \c^n$ be a smooth quasi projective
variety and $L(x_1,\dots,x_n)$ a general linear form.
Then $L|_X$ has only nondegenerate critical points. \qed
\end{lem}

\begin{say}{\ }
It is sufficient to  prove (\ref{real.loj.max.thm})
after a linear change of coordinates and after
replacing the $f_i$ by $\epsilon_if_i$. Thus, in order to simplify
notation, we may assume that for every $J\subset I$, 
 $0<|\eta|<\eta_0$ and $0<\delta\ll 1$ the following are
satisfied:
\begin{enumerate}
\item $X_{J,\eta}:=\{{\bold x}\in \c^n|f_j({\bold x})=\eta\ \forall
j\in J\}$ is a smooth complete intersection of dimension $n-|J|$.
\item  The infimum of $\Phi$ on the real cube $\max_i|x_i|=\delta$
is achieved at an interior point of the face $x_1=\delta$.
\item  $x_1$ has only nondegenerate critical points
on $X_{J,\eta}$.
\end{enumerate}
Set
$D^*_J:=\{{\bold x}| f_j({\bold x})=f_{j'}({\bold x})\ \forall
j,j'\in J\qtq{and}
\rank (\jac_J)\leq |J|\}$.
$D^*_J$ can have irreducible components of  dimension bigger than
one, so let $C^*_J\subset D^*_J$ denote the union of all those
irreducible components which intersect $X_{J,\eta}$ for some
$0<|\eta|<\eta_0$. Then $C^*_J$ is an algebraic curve and we
want to prove an estimate
$$
\Phi({\bold x})\geq C\cdot ||{\bold x}||^{N}
\qtq{for ${\bold x}\in C^*_J$}.
$$ 
The following lemma reduces the computation of $N$
to computing the number of intersection points
of $ C^*_J$ with the hypersurface $(f_j=\eta)$.
Since $ C^*_J\subset \cup_{\eta}X_{J,\eta}$, 
these intersections are precisely the critical points
of $x_1$ on $X_{J,\eta}$, hence they do  not depend on
the choice of $j\in J$.

It is worth noting that the number of intersection points
may be less than $d\cdot \deg C^*_J$ since there can be
intersection points at infinity.
\end{say}

\begin{lem} Let $D\subset \c^n$ be an irreducible   algebraic
or analytic 
curve and $f$   a
polynomial  or analytic function   not identically zero on $D$
such that $f(0)=0$. Let $N=N_{\epsilon}$ denote the number of 
solutions of $f({\bold x})=\epsilon, {\bold x}\in D$
for general $|\epsilon|\ll 1$. 

Then $|f({\bold x})|\geq C\cdot ||{\bold x}||^{N}$
for ${\bold x}\in D$ near $0$.
\end{lem}

Proof. Let $\pi:\bar D\to D$ be the normalization 
and $t$ a local
parameter on $\bar D$ at a point  $q\in \pi^{-1}(0)$.
$\pi$ is given by its analytic coordinate functions $\pi_j(t)$.
Let $M$ be the order of vanishing  of $f\circ \pi$ at $q$.
Then $|f(\pi(t))|\geq C\cdot |t|^M$
and $||\pi(t)||\sim |t|^m$ where $m=\min_i\{\mult_q\pi_i(t)\}$.
Thus $|f(\pi(t))|\geq C\cdot ||\pi(t)||^{M/m}$.
On the other hand, $M\leq N$.\qed
\medskip

The next result gives a formula for the number of critical points
of a  function on a 
 smooth complete 
intersection.

\begin{prop}\label{crit.pt.formula}
 Let $X\subset \a^n$ be a smooth complete 
intersection of hypersurfaces of degrees $d_1,\dots,d_k$.
Let $g$ be a polynomial of degree $c$ with only isolated
critical points on $X$. Then the number of critical points
is at most the coefficient of $H^n$
in the Taylor series of 
$$
(-1)^{n-k} \frac{(1+H)^n}{1+cH}\prod_i
\frac{d_iH}{1+d_iH}.
$$
Equality holds for general choice of the hypersurfaces and $g$.
\end{prop}

Proof. Let $Y_t:=(g=t)\subset X$ be a  level set.
Let $\bar X\subset \p^n$ denote the closure of $X$ and set
$X_{\infty}:=\bar X\setminus X$. Define $\bar Y_t$ and
$Y_{t,\infty}$ analogously. Let  $x_0=0$ be the equation of the
hyperplane at infinity. 
$g:\bar X\map \p^1$ is a rational map
which becomes a morphism 
$\tilde g:\tilde X:=B_Z\bar X\to \p^1$
after blowing up the base locus $Z:=(g=x_0^c=0)$.
(This needs to be done scheme theoretically.)
$\tilde g^{-1}(t)\cong \bar Y_t$ for $t\in \a^1$.
We can identify $X$ with a subset of $\tilde X$,
set $\tilde X_{\infty}:=\tilde X\setminus X$.
Under some mild conditions the homology of a pair
$(A,B)$ depends only on $A\setminus B$ (cf.
\cite[19.14]{green-harp}). Thus (\ref{euler.char.lem}) gives that
$$
\begin{array}{l}
e(\tilde X,\tilde X_{\infty})=e(\bar X,X_{\infty})=\\
e(\bar Y,Y_{\infty})e(\c\p^1,\{\infty\})
+\sum_{t\in \c}[e(\bar Y_t,Y_{t,\infty})-e(\bar Y,Y_{\infty})],
\end{array}
$$
where $Y$ denotes a typical level set $Y_t$. 
If $\tilde g$ has a single nondegenerate critical point 
on $\bar Y_t$ then
by the Picard--Lefschetz formulas (see, for instance,
\cite[Chap.3]{loo})
$$
e(\bar Y_t,Y_{t,\infty})-e(\bar Y,Y_{\infty})=(-1)^{\dim Y+1}.
$$
Thus if $\tilde g$
has only nondegenerate critical points over $\a^1$
then their
  number is 
$$
\begin{array}{l}
(-1)^{n-k}[e(\bar X,X_{\infty})-e(\bar Y,Y_{\infty})]=\\
 (-1)^{n-k}[e(\bar X)-e(X_{\infty})-e(\bar Y)+e(Y_{\infty})].
\end{array}
$$
Assume that $\bar X,X_{\infty}, \bar Y, Y_{\infty}$ are all
smooth complete intersections. The total Chern class of a complete
intersection of type $b_1,\dots,b_m$ in $\p^n$ is given by
$(1+H)^{n+1}/\prod (1+b_iH)$ (cf.\ \cite[\S 22]{hirz}).
Substituting these into the above formula we obtain the result.

For arbitrary $X$ and $g$ we argue that the number of
critical points (in the complex domain)
 is an upper semi continuous functions
(provided that it is finite). This follows from the 
observation that if $g_s$ is a family of smooth functions
depending smoothly on a parameter $s$ and $g_0$ has an isolated
critical point at   $p_0$ then each $g_s$
has at least one isolated
critical point near $p_0$. Indeed, an isolated critical point is
an isolated solution of the equations $\partial g_0/\partial
x_j=0\ \forall j$. These are $\dim X$ equations in $\dim X$
variables,  thus the system
$\partial g_t/\partial x_j=0\ \forall j$ has a nearby solution
for small
$|t|$.
\qed

\begin{cor}\label{expl.formula}
 Let $X\subset \a^n$ be a smooth complete 
intersection of $k$ hypersurfaces of degree $d$.
Let $g$ be a linear polynomial  with only isolated
critical points on $X$. Then the number of critical points
is at most
$$
\binom{n-1}{k-1}d^k(d-1)^{n-k}.
$$
Equality holds for general choice of the hypersurfaces and $g$.
\end{cor}

Proof. By (\ref{crit.pt.formula}) our upper bound
is the coefficient of $H^n$ in 
$$
\begin{array}{l}
(-1)^{n-k}(1+H)^{n-1}
\left(\frac{dH}{1+dH}\right)^k=\\
(-1)^{n-k}d^k
H^k\left(\sum \binom{n-1}{s}H^s\right)
\left(\sum (-1)^sd^s\binom{s+k-1}{k-1}H^s\right).
\end{array}
$$
Thus the coefficient of $H^n$ is
$$
\begin{array}{l}
(-1)^{n-k}d^k\sum_{s=0}^{n-k}\binom{n-1}{s}
(-d)^{n-k-s}\binom{n-s-1}{k-1}=\\
d^n\sum_{s=0}^{n-k}\binom{n-1}{s}\binom{n-s-1}{k-1}
\left(\frac{-1}{d}\right)^s=\\
d^n\binom{n-1}{k-1}\sum_{s=0}^{n-k}\binom{n-k}{s}
\left(\frac{-1}{d}\right)^s=\\
\binom{n-1}{k-1}d^n\left(1-\frac{1}{d}\right)^{n-k}.\qed
\end{array}
$$

The next result is  known in various forms, see, for instance,
\cite[III.11.4]{bpv}.

\begin{lem}\label{euler.char.lem} Let $A$ be a compact simplicial
complex, $B\subset A$ a subcomplex, $C$ a  topological surface
(possibly with boundary) and
$D\subset C$ a finite set.
 Let
$g:A\to C$ be a    continuous map such that $g^{-1}(D)\subset
B$. For $c\in C$ set $A_c:=g^{-1}(c)$ and $B_c:=g^{-1}(c)\cap B$.
Assume that
$g$ is a locally trivial fibration of the pair $(A,B)$
with typical fiber $(A_{gen},B_{gen})$
except over  finitely many points.
Then
$$
e(A,B)=e(C,D)e(A_{gen},B_{gen})+\sum_{c\in C\setminus D}
[e(A_c,B_c)-e(A_{gen},B_{gen})].
$$
\end{lem}

Proof. If $(A,B)\to Z$ is a locally trivial fibration 
with typical fiber $(F,E)$ then $e(A,B)=e(Z)e(F,E)$.
Thus any locally trivial fibration  over $S^1$ has
zero Euler number. The Mayer--Vietoris sequence now shows
that if $(A,B)=(A_1,B_1)\cup (A_2,B_2)$
and $(A_1,B_1)\cap (A_2,B_2)$ is a 
locally trivial fibration  over $S^1$ then
$e(A,B)=e(A_1,B_1)+ e(A_2,B_2)$. By repeatedly cutting discs out
of $C$, the proof of (\ref{euler.char.lem}) is  reduced to the
following three special cases:
\begin{enumerate}
\item $D=\emptyset$ and $(A,B)\to C$ is a locally trivial
fibration. In this case the formula holds, as we noted above.
\item $C$ is a disc, $D$ is empty 
and  $(A,B)\to C$ is a locally trivial
fibration over $C\setminus p$. Then  $(A,B)$ retracts
to  $(A_p,B_p)$ and so $e(A,B)=e(A_p,B_p)$.\qed
\item $C$ is a disc, $D$ is a single point $p$
and  $(A,B)\to C$ is  locally trivial
 over $C\setminus p$. Then  
  $A_p=B_p$ and so $e(A,B)=e(A_p,B_p)=0$.\qed
\end{enumerate}

\begin{ack}  I   thank B. Fulton, J.\ Gwo\'zdziewicz, T. Krasi\'nski,
 T.\ Kuwata and R.\ Pemantle
for helpful
comments and references. J.\ Johnson  helped me  with
Maple computations relating to these results.
Partial financial support was provided by  the NSF under grant number 
DMS-9622394. 
\end{ack}

\vskip1cm

\noindent University of Utah, Salt Lake City UT 84112 

\begin{verbatim}kollar@math.utah.edu\end{verbatim}


\begin{thebibliography}{Chislenko88}




\bibitem[BPV84]{bpv}  W.   Barth,  C.   Peters and A.   Van de Ven,
  Compact Complex Surfaces, Springer 1984 

\bibitem[Brownawell88]{brownaw-loc} W. D. Brownawell, Local
diophantine  Nullstellen inequalities, Journal AMS  1 (1988)
 311-322



\bibitem[CKT99]{ckt} E. Cygan, T. Krasi\'nski and P.
Tworzewski, Separation at infinity and the \L ojasiewicz
exponent of polynomial mappings, Invent Math. 
136 (1999) 75--88


\bibitem[Fulton84]{fulton}  W.   Fulton, Intersection Theory,
 Springer  (1984)

\bibitem[Greenberg-Harper81]{green-harp} M. Greenberg and
J. Harper, Algebraic topology, Benjamin/Cummings, 1981 

\bibitem[Gwo\'zdziewicz99a]{gwoz1} J. Gwo\'zdziewicz, The 
\L ojasiewicz exponent of an analytic function at an isolated
zero, Comment. Math. Helv. (to appear)

\bibitem[Gwo\'zdziewicz99b]{gwoz2} J. Gwo\'zdziewicz, Growth at
infinity of a polynomial with a compact zero set, 
Singularities symposium --  \L ojasiewicz 70, 
Banach Cent.\ Publ. Warszawa, 1998

\bibitem[HLP34]{hlp} G. Hardy, J. Littlewood and G. P\'olya,
Inequalities, Cambridge Univ. Press, 1934


\bibitem[Hartshorne77]{harts-book}   R.   Hartshorne,
Algebraic Geometry, Springer,  1977

\bibitem[Hirzebruch66]{hirz}   F. Hirzebruch, Topological
methods in algebraic geometry, third enlarged ed.,  Springer,    1966



\bibitem[JKS92]{jks}  S. Ji,   J. Koll\'ar and  B.  Shiffman,  A
Global \L ojasiewicz Inequality for Algebraic Varieties, 
Trans. AMS,  329 (1992) 813-818 

\bibitem[Koll\'ar88]{null} J.   Koll\'ar, Sharp effective
Nullstellensatz,  Jour. AMS, 1 (1988)
 963-975

\bibitem[\L ojasiewicz59]{loj-stud}   S. \L ojasiewicz,  Sur
le probl\`eme de la division,
 Studia Math  18 (1959)  87-136 


\bibitem[Looijenga84]{loo} E. Looijenga, Isolated singular points
of complete intersections, Cambridge Univ.\ Press, 1984

\bibitem[Looijenga97]{loo-pc} E. Looijenga, Cohomology and intersection
homology of algebraic varieties, in: Complex algebraic geometry, 
IAS/Park City math.\ series vol.\ 3. AMS,  1997, 221--264

\bibitem[Miyaoka84]{miya}  Y.\ Miyaoka, The maximal number of
 quotient singularities on
surfaces with given numerical invariants. Math. Ann. 268 (1984)
159--171

\bibitem[Solern\'o91]{solerno} P. Solern\'o,   Effective
 \L ojasiewicz inequalities in
semialgebraic geometry. Appl. Algebra Engrg. Comm. Comput. 2 (1991)
2--14.

\bibitem[Varchenko83]{var} 
 A. N. Varchenko, Semicontinuity of the spectrum and an upper
bound for the number of singular points of the projective hypersurface.
 (Russian)
Dokl. Akad. Nauk SSSR 270 (1983) 1294--1297

\bibitem[Yoshihara79]{yoshi}  H.\ Yoshihara, On Plane Rational Curves,
 Proc. Japan Acad. 55  (1979) 152--155

\end{thebibliography}
\end{document}